\documentclass[a4paper,11pt]{article}

\usepackage{latexsym}
\usepackage{amssymb,amsmath}
\usepackage{amsthm}
\usepackage{xcolor}
\usepackage{ulem}
\usepackage{graphicx}
\usepackage{tikz} %for tikzｗ
\usetikzlibrary{positioning, intersections, calc, arrows.meta,math} %tikz library
\newtheorem{theorem}{Theorem}
\newtheorem{corollary}[theorem]{Corollary}
\newtheorem{conjecture}[theorem]{Conjecture}

\newtheorem{lemma}[theorem]{Lemma}

\newcommand{\pr}{\noindent{\bf Proof. \hspace{0.5em}}}

\date{\empty}

\begin{document}

%%%%%%%%%%%%%%%%%%%%%%%%%%%%%%%%%%%%%%%%%%%%%%%%%%%%%%%%%%%%%%%
%%%%%%%%%%%%%%%%%%%%%%%%%%%%%%%%%%%%%%%%%%%%%%%%%%%%%%%%%%%%%%%
\title{
A generalization of an ear decomposition and
 $k$-trees in highly connected star-free graphs
}
%%%%%%%%%%%%%%%%%%%%%%%%%%%%%%%%%%%%%%%%%%%%%%%%%%%%%%%%%%%%%%%
%%%%%%%%%%%%%%%%%%%%%%%%%%%%%%%%%%%%%%%%%%%%%%%%%%%%%%%%%%%%%%%

\author{
Shun-ichi Maezawa$^{1}$
\thanks{This work was partially supported by JSPS KAKENHI Grant Number JP22K13956}
\thanks{E-mail address: \texttt{maezawa.mw@gmail.com}}
\hspace{+8pt}
Kenta Ozeki$^{2}$
\thanks{This work was partially supported by JSPS KAKENHI Grant Numbers JP22K19773 and JP23K03195 }
\thanks{E-mail address: \texttt{ozeki-kenta-xr@ynu.ac.jp}}
\hspace{+8pt}
 \\
Masaki Yamamoto$^{3}$
\thanks{E-mail address: \texttt{yamamoto@st.seikei.ac.jp}}
\hspace{+8pt}
Takamasa Yashima$^{3,4}$
\thanks{This work was partially supported by JSPS KAKENHI Grant Number JP20K14353}
\thanks{E-mail address: \texttt{takamasa.yashima@gmail.com}}
\vspace{+8pt}
 \\
\small%%%%%
$^1$\small\textsl{
Department of Information Science, Nihon University,}\\
\small\textsl{Sakurajosui 3-25-40, Setagaya-Ku, Tokyo, 156-8550, Japan}
\vspace{+6pt}\\
\small%%%%%
$^2$\small\textsl{
Faculty of Environment and Information Sciences, Yokohama National University,}\\
\small\textsl{79-2 Tokiwadai, Hodogaya-ku, Yokohama, Kanagawa, 240-8501, Japan}
\vspace{+6pt}\\
\small%%%%%
$^3$\small\textsl{
Department of Computer and Information Science, Seikei University,}\\
\small\textsl{3-3-1 Kichijoji-Kitamachi, Musashino-shi, Tokyo, 180-8633, Japan}
\vspace{+6pt}\\
\small%%%%%
$^4$\small\textsl{
Kanazawa Institute of Technology,}\\
\small\textsl{7-1 Ohgigaoka, Nonoichi-shi, Ishikawa, 921-8501, Japan}
}

\maketitle

\vspace{-1cm}

%\begin{center} 
%\hspace{2em} \date{\today}
%\end{center}

%%%%%%%%%%%%%%%%%%%%%%%%%%%%%%%%%%%%%%%%%%%%%%%%%%%%%%%%%%%%%%%
\begin{abstract}
In this paper, we introduce a generalized version of an ear decomposition,
called a $j$-spider decomposition, for $j$-connected star-free graphs with~$j \geq 2$.
Its application enables us to
improve a previousely known sufficient condition for the existence of a~$k$-tree
in highly connected star-free graphs,
where a~$k$-tree is a spanning tree in which every vertex is of degree at most~$k$. 
More precisely, we show that every $j$-connected $K_{1,j(k-2)+2}$-free graph has a~$k$-tree for~$k\ge j$,
thereby improving a classical result of Jackson and Wormald for~$k\ge j$.
Our approach differs from previous studies based on toughness-type arguments
and instead relies on both a~$j$-spider decomposition and
a factor theorem related to Hall's marriage theorem.
\end{abstract}
%%%%%%%%%%%%%%%%%%%%%%%%%%%%%%%%%%%%%%%%%%%%%%%%%%%%%%%%%%%%%%%



%%%%%%%%%%%%%%%%%%%%%%%%%%%%%%%%%%%%%%%%%%%%%%%%%%%%%%%%%%%%%%%
%%%%%%%%%%%%%%%%%%%%%%%%%%%%%%%%%%%%%%%%%%%%%%%%%%%%%%%%%%%%%%%
\section{Introduction}
%%%%%%%%%%%%%%%%%%%%%%%%%%%%%%%%%%%%%%%%%%%%%%%%%%%%%%%%%%%%%%%
%%%%%%%%%%%%%%%%%%%%%%%%%%%%%%%%%%%%%%%%%%%%%%%%%%%%%%%%%%%%%%%

Throughout this paper, %note,
we consider only finite, simple, and undirected graphs.
Let~$G$ be a graph.
We let~$V(G)$ and~$E(G)$ denote the vertex set and the edge set of~$G$, respectively.
For~$x \in V (G)$, $d_G(x)$ denotes the {\it degree} of~$x$ in~$G$.

An {\it ear decomposition},
which is defined below,
is a well-known and powerful method for studying $2$-connected graphs
(e.g.,~for matching theory~\cite{L2} and \cite[Chapters 4 and 9.1]{LP},
a strong orientation \cite[Theorem 5.10]{BM},
and Hamiltonian cycles in square graphs \cite{AGRT}).
However, less is known about corresponding techniques for graphs with higher connectivity.
In this paper, we first present a use of an ear decomposition to star-free graphs,
propose its generalization to $j$-connected star-free graphs with~$j \geq 2$,
and then demonstrate its application to $k$-trees.

%%%%%%%%%%%%%%%%%%%%%%%%%%%%%%%%%%%%%%%%%%%%%%%%%%%%%%%%%%%%%%%
\subsection{Ear decomposition of star-free graphs}
%%%%%%%%%%%%%%%%%%%%%%%%%%%%%%%%%%%%%%%%%%%%%%%%%%%%%%%%%%%%%%%

A {\it non-trivial ear} of a graph~$G$ with respect to a vertex subset~$X$ 
is a path of length at least~$2$ connecting two vertices in~$X$
such that no internal vertices belong to~$X$.
(Such a path without the length condition is called an {\it ear}.
However, we consider only non-trivial ones in this paper,
since we require that all vertices belong to at least one ear,
but do not mind edges belonging to no ears.)
It is well-known that for any vertex subset~$X$ with~$|X|\ge 2$ of a $2$-connected graph~$G$
and any vertex~$v \in V(G) \setminus X$,
there exists a non-trivial ear with respect to~$X$ passing through~$v$ in which the end vertices are distinct.
This property enables us to give an {\it ear decomposition} of~$G$,
that is,
a sequence of subgraphs $C_0, C_1, \dots, C_m$ such that
$C_0$ is a cycle of~$G$,
each $C_{\ell}$ with $1 \leq \ell \leq m$
is a non-trivial ear with respect to $\bigcup_{i=0}^{\ell-1} V(C_{i})$,
and $V(G) = \bigcup_{i=0}^{m} V(C_{i})$.

We say that a graph $G$ is {\it $K_{1,t}$-free} or {\it $t$-star-free}
if $G$ contains no $K_{1,t}$ as an induced subgraph,
where $K_{1,t}$ is the complete bipartite graph with partite sets of cardinalities~$1$ and~$t$.
For $2$-connected $K_{1,t}$-free graphs,
Ku\v{z}el and Teska~\cite{KT} proved the following theorem.

%%%%%%%%%%%%%%%%%%%%%%%%%%%%%%%%%%%%%%%%%%%%%%%%%%%%%%%%%%%%%%%
\begin{theorem}[Ku\v{z}el and Teska~\cite{KT}]
\label{KTthm}
Let $t\ge 2$ be an integer.
Let~$G$ be a~$2$-connected $K_{1,t}$-free graph.
Then $G$ has a~$2$-connected spanning subgraph with maximum degree at most~$t$.
\end{theorem}
%%%%%%%%%%%%%%%%%%%%%%%%%%%%%%%%%%%%%%%%%%%%%%%%%%%%%%%%%%%%%%%

Note that their proof is based on an analysis of $2$-connected graphs,
and does not make use of an ear decomposition.
We first give an alternative proof in Section~\ref{proofj=2_sec}
by taking an appropriate ear decomposition.
While the proof in~\cite{KT} spends five pages,
our proof requires only a half page.
In addition, we will extend this idea to prove our main theorem stated in the next subsection.

%%%%%%%%%%%%%%%%%%%%%%%%%%%%%%%%%%%%%%%%%%%%%%%%%%%%%%%%%%%%%%%
\subsection{Generalized version of an ear decomposition for highly connected graphs}
%%%%%%%%%%%%%%%%%%%%%%%%%%%%%%%%%%%%%%%%%%%%%%%%%%%%%%%%%%%%%%%

As a natural extension of a non-trivial ear, we define a~{\it $j$-spider} as follows:
A tree having at most one vertex of degree greater than two is called a~{\it spider}.
The vertex of degree greater than two is called its {\it branch} if exists;
otherwise, any vertex of degree two can be regarded as its {\it branch}.
Let $j \geq 2$ be an integer.
In a spider, 
its leaf is particularly called a~{\it foot},
and a spider with exactly~$j$ feet is called a~{\it $j$-spider}.
In other words,
a spider is a~$j$-spider if and only if the branch has degree exactly~$j$.
It can also be defined as a subdivision of~$K_{1,j}$.
For a vertex~$v$ and a vertex subset~$X$ of a graph~$G$ with~$v \notin X$ and $|X|\ge j$,
a~$j$-spider {\it from $v$ to~$X$}
is one in which
$v$ is its branch,
all feet are contained in~$X$ with any pair of the feet distinct and 
no other vertices belong to~$X$.
By Menger's theorem,
for any vertex subset~$X$ of a~$j$-connected $G$ with~$|X| \geq j$
and any vertex~$v \notin X$,
there exists a~$j$-spider from~$v$ to~$X$.

Similarly to an ear decomposition of a~$2$-connected graph,
in some sense,
we can take a~``$j$-spider decomposition'' of a~$j$-connected graph.
As in the following theorem and its application explained later,
it is a powerful method
to find a spanning subgraph with particular degree conditions
in~$j$-connected star-free graphs.

%%%%%%%%%%%%%%%%%%%%%%%%%%%%%%%%%%%%%%%%%%%%%%%%%%%%%%%%%%%%%%%
\begin{theorem}
\label{spider_mainthm}
Let $j\ge 2$ and $t\ge2$ be integers.
Let $G$ be a $j$-connected $K_{1,t}$-free graph.
Then $G$ has a sequence of subgraphs $F_0, F_1, \dots, F_m$
such that $F_0$ is a cycle with~$|V(F_0)| \geq j$,
each $F_{\ell}$ with~$1 \leq \ell \leq m$
is a~$j$-spider from a vertex in $G - \bigcup_{i=0}^{\ell-1} V(F_{i})$ to $\bigcup_{i=0}^{\ell-1} V(F_{i})$,
$V(G)=\bigcup_{i=0}^{m}V(F_i)$,
and every vertex $v$ satisfies the following property:
\begin{itemize}
\item
If $j \geq 3$ and there exists a $j$-spider $F_{\ell}$ (with $1 \leq \ell \leq m$)
whose branch is $v$,
then $v$ is a foot of none of the $j$-spiders $F_0, F_1, \dots, F_m$.
\item
If $j=2$ or 
there exists no $j$-spider $F_{\ell}$ (with $1 \leq \ell \leq m$)
whose branch is $v$,
then $v$ is a foot of at most $t-2$ of the $j$-spiders $F_0, F_1, \dots, F_m$.
\end{itemize}
\end{theorem}
%%%%%%%%%%%%%%%%%%%%%%%%%%%%%%%%%%%%%%%%%%%%%%%%%%%%%%%%%%%%%%%

In addition,
by taking the union of~$F_0, F_1, \dots, F_m$,
Theorem~\ref{spider_mainthm} immediately yields the following consequence.

%%%%%%%%%%%%%%%%%%%%%%%%%%%%%%%%%%%%%%%%%%%%%%%%%%%%%%%%%%%%%%%
\begin{corollary}
%[Ku\v{z}el and Teska~\cite{KT}]
\label{KTthm+}
Let $j\ge 2$ and $t\ge2$ be integers.
Let $G$ be a $j$-connected $K_{1,t}$-free graph.
Then $G$ has a $2$-connected spanning subgraph
with maximum degree at most~$\max\{j,t\}$.
\end{corollary}
%%%%%%%%%%%%%%%%%%%%%%%%%%%%%%%%%%%%%%%%%%%%%%%%%%%%%%%%%%%%%%%

Note that the case~$j=2$ of Corollary \ref{KTthm+} is equivalent to Theorem \ref{KTthm}.
Since our proof of the case~$j=2$ of Theorem~\ref{spider_mainthm} is based on the same idea as,
but simpler than, the general case~$j \geq 3$,
we will give it separately in Subsection~\ref{proofj=2_sec},
which also gives an alternative proof to Theorem \ref{KTthm}.

%%%%%%%%%%%%%%%%%%%%%%%%%%%%%%%%%%%%%%%%%%%%%%%%%%%%%%%%%%%%%%%
\subsection{Application to the existence of a $k$-tree}
%%%%%%%%%%%%%%%%%%%%%%%%%%%%%%%%%%%%%%%%%%%%%%%%%%%%%%%%%%%%%%%

We now describe an application of Theorem~\ref{spider_mainthm}
to~$k$-trees.
Let $k\ge2$ be an integer, and~$G$ be a connected graph.
A {\it $k$-tree} of~$G$ is defined as a spanning tree~$T$ such that $d_T(v)\le k$ for all $v\in V(T)$.
A {\it $k$-walk} is a closed spanning walk visiting each vertex at most~$k$ times.
Note that a~$2$-tree and a~$1$-walk correspond to a Hamiltonian path and a Hamiltonian cycle, respectively.

In 1990, Jackson and Wormald~\cite{JW} 
showed that
for~$k \ge 2$,
if a graph~$G$ has a~$k$-tree,
then $G$ has a~$k$-walk,
and if $G$ has a~$(k-1)$-walk,
then $G$ has a~$k$-tree.
Moreover, they 
gave a sufficient condition
for star-free graphs to have a~$k$-tree and a~$k$-walk,
and showed that the star-free condition 
can be relaxed with respect to the connectivity of graphs under consideration.

%%%%%%%%%%%%%%%%%%%%%%%%%%%%%%%%%%%%%%%%%%%%%%%%%%%%%%%%%%%%%%%
\begin{theorem}[Jackson and Wormald~\cite{JW}]
\label{JWthm2}
Let $k\ge 2$ be an integer.
Let $G$ be a connected $K_{1,k}$-free graph.
Then $G$ has a~$(k-1)$-walk (and hence a~$k$-tree).
\end{theorem}
%%%%%%%%%%%%%%%%%%%%%%%%%%%%%%%%%%%%%%%%%%%%%%%%%%%%%%%%%%%%%%%

%%%%%%%%%%%%%%%%%%%%%%%%%%%%%%%%%%%%%%%%%%%%%%%%%%%%%%%%%%%%%%%
\begin{theorem}[Jackson and Wormald~\cite{JW}]
\label{JWthm3}
Let $j\ge 1$ and $k\ge 3$ be integers.
Let $G$ be a~$j$-connected $K_{1,j(k-2)+1}$-free graph.
Then $G$ has a~$k$-tree
(and hence a~$k$-walk).
\end{theorem}
%%%%%%%%%%%%%%%%%%%%%%%%%%%%%%%%%%%%%%%%%%%%%%%%%%%%%%%%%%%%%%%

Note that Theorem~\ref{JWthm2} for~$k$-trees was already shown in~\cite{CKR},
and that the case~$j=1$ for Theorem~\ref{JWthm3} is weaker than Theorem~\ref{JWthm2}.
By using Theorem~\ref{spider_mainthm},
we show the following result,
which 
generalizes Theorem~\ref{JWthm3} for $k\ge j$.

\color{black}
%%%%%%%%%%%%%%%%%%%%%%%%%%%%%%%%%%%%%%%%%%%%%%%%%%%%%%%%%%%%%%%
\begin{theorem}
\label{main}
Let $j\ge 1$ and $k\ge j$ be integers.
Let $G$ be a~$j$-connected $K_{1,j(k-2)+2}$-free graph.
Then $G$ has a~$k$-tree.
\end{theorem}
%%%%%%%%%%%%%%%%%%%%%%%%%%%%%%%%%%%%%%%%%%%%%%%%%%%%%%%%%%%%%%%

In previous studies, the existence of $k$-trees in star-free graphs has typically been established in terms of {\it toughness}.
The {\it toughness} of a connected non-complete graph~$G$, denoted by $t(G)$, is defined as
$t(G):=\min\{  \frac{|X|}{\omega(G-X)} \mid \emptyset\neq X\subseteq V(G), \omega(G-X) \ge 2 \}$,
where $\omega(G-X)$ denotes the number of components of~$G-X$.
A graph~$G$ is said to be {\it $t$-tough} if $t(G) \ge t$.
In fact, the proof of Theorem~\ref{JWthm3} in~\cite{JW} was done by
showing that 
any~$j$-connected $K_{1,j(k-2)+1}$-free graph
is $\frac{1}{k-2}$-tough,
and hence there exists a~$k$-tree due to Win's toughness-type theorem~\cite{Win},
where $j\ge 1$ and $k\ge 3$.
Note that there are $j$-connected $K_{1,j(k-2)+2}$-free graphs that are not 
$\frac{1}{k-2}$-tough (e.g.,~the complete graph $K_{j, j(k-2)+1}$),
and hence such an approach does not seem applicable to the proof of Theorem~\ref{main} %our main theorem.

On the other hand, as we wrote above,
the key idea to improve Theorem~\ref{JWthm3}
is Theorem~\ref{spider_mainthm},
together with a factor theorem related to Hall's marriage theorem~(see Section~\ref{proofmain_sec}),
which is different from toughness-type arguments mentioned above.

Jackson and Wormald posed the 
following conjecture concerning $k$-walks.

%%%%%%%%%%%%%%%%%%%%%%%%%%%%%%%%%%%%%%%%%%%%%%%%%%%%%%%%%%%%%%%
\begin{conjecture}[Jackson and Wormald~\cite{JW}]
\label{JWconj1}
Let $j\ge 1$ and $k\ge 2$ be integers.
Let $G$ be a~$j$-connected $K_{1,j(k-1)+1}$-free graph.
Then $G$ has a~$(k-1)$-walk.
\end{conjecture}
%%%%%%%%%%%%%%%%%%%%%%%%%%%%%%%%%%%%%%%%%%%%%%%%%%%%%%%%%%%%%%%

However, in 2004,
Jin and Li~\cite{JL} provided counterexamples showing that Conjecture~\ref{JWconj1} is false for~$j\ge 2$,
while it clearly holds for~$j=1$ by Theorem~\ref{JWthm2}.
If Conjecture~\ref{JWconj1} were true,
then 
we would show the existence of a~$k$-tree in any graph satisfying the assumption of Conjecture~\ref{JWconj1},
because any graph containing a~$(k-1)$-walk has a~$k$-tree.
Thus, we can formulate the following conjecture concerning $k$-trees~(noting that
the above examples, due to Jin and Li,
are $j$-connected $K_{1,j(k-1)+1}$-free graphs that contain no $(k-1)$-walks but do contain $k$-trees).
Theorem~\ref{main} is a partial affirmative solution to Conjecture~\ref{JWconj2}.

%%%%%%%%%%%%%%%%%%%%%%%%%%%%%%%%%%%%%%%%%%%%%%%%%%%%%%%%%%%%%%%
\begin{conjecture}%[Jackson and Wormald \cite{JW}]
\label{JWconj2}
Let $j\ge 1$ and $k\ge 2$ be integers.
Let $G$ be a $j$-connected $K_{1,j(k-1)+1}$-free graph.
Then $G$ has a~$k$-tree.
\end{conjecture}
%%%%%%%%%%%%%%%%%%%%%%%%%%%%%%%%%%%%%%%%%%%%%%%%%%%%%%%%%%%%%%%

%%%%%%%%%%%%%%%%%%%%%%%%%%%%%%%%%%%%%%%%%%%%%%%%%%%%%%%%%%%%%%%
\subsection{Further application to spanning trees with degree constraint}
\label{introtree_sec}
%%%%%%%%%%%%%%%%%%%%%%%%%%%%%%%%%%%%%%%%%%%%%%%%%%%%%%%%%%%%%%%

The idea of the proof of Theorem~\ref{main}
can be used to 
a spanning tree with bounded degrees.
The following theorem was known.

%%%%%%%%%%%%%%%%%%%%%%%%%%%%%%%%%%%%%%%%%%%%%%%%%%%%%%%%%%%%%%%
\begin{theorem}[Hasanvand~\cite{Hasanvand} (see also~\cite{CS, ENV, LX})]
\label{maxdeg_thm}
Let $j\ge 1$ be an integer.
Let $G$ be a~$j$-edge-connected graph.
Then $G$ has a spanning tree~$T$ 
with $d_T(v) \leq \lceil \frac{d_G(v) -2}{j}\rceil +2$ for any vertex~$v$.
\end{theorem}
%%%%%%%%%%%%%%%%%%%%%%%%%%%%%%%%%%%%%%%%%%%%%%%%%%%%%%%%%%%%%%%

Note that Theorem~\ref{maxdeg_thm} and the weaker versions were shown
by using toughness-type arguments in~\cite{CS, ENV, Hasanvand}
and a more naive approach in~\cite{LX}.
We point out that
by employing the same arguments as in the proof of Theorem~\ref{main},
we can show the existence of a spanning tree~$T$
in a~$j$-connected graph
such that $d_T(v) \leq \max\left\{\lceil \frac{d_G(v) -2}{j}\rceil +2, j\right\}$
for any vertex~$v$.
Although we need the~$j$-connected assumption
instead of the~$j$-edge-connected one,
our proof is based on Theorem \ref{spider_mainthm},
which is totally different from those in~\cite{CS, ENV, Hasanvand, LX}.
We will explain more after the proof of Theorem~\ref{main}.

%%%%%%%%%%%%%%%%%%%%%%%%%%%%%%%%%%%%%%%%%%%%%%%%%%%%%%%%%%%%%%%
%%%%%%%%%%%%%%%%%%%%%%%%%%%%%%%%%%%%%%%%%%%%%%%%%%%%%%%%%%%%%%%
\section{Preliminaries}
%%%%%%%%%%%%%%%%%%%%%%%%%%%%%%%%%%%%%%%%%%%%%%%%%%%%%%%%%%%%%%%
%%%%%%%%%%%%%%%%%%%%%%%%%%%%%%%%%%%%%%%%%%%%%%%%%%%%%%%%%%%%%%%

In this section,
we prepare the necessary terminology and
state a key lemma
that will be used in the proof of Theorem~\ref{spider_mainthm}.

Our notation is standard, and is mostly taken from Diestel~\cite{D25}.
Possible exceptions are as follows:
Let~$G$ be a graph.
For~$x \in V(G)$, $N_G(x)$ denotes the set of vertices adjacent to~$x$ in~$G$;
thus $d_G(x)=|N_G(x)|$.
For $S \subseteq V(G)$,
the graph obtained from~$G$ by deleting all vertices in $S$ together with the edges incident with them is denoted by $G-S$.
A vertex~$x$ of~$G$ is often identified with the set~$\{x\}$;
for example, if $H$ is a subgraph with~$x\in V(H)$,
then we write $H-x$ for~$H-\{x\}$.

Using $2$-spiders (i.e., non-trivial ears),
we first show the case~$j=2$ of Theorem~\ref{spider_mainthm}
in the next subsection.
Since it immediately yields Theorem~\ref{KTthm},
this gives its alternative proof,
which is shorter than the one in~\cite{KT}.

%%%%%%%%%%%%%%%%%%%%%%%%%%%%%%%%%%%%%%%%%%%%%%%%%%%%%%%%%%%%%%%
\subsection{
Proof of the case~$j=2$ of Theorem~\ref{spider_mainthm}}
\label{proofj=2_sec}
%%%%%%%%%%%%%%%%%%%%%%%%%%%%%%%%%%%%%%%%%%%%%%%%%%%%%%%%%%%%%%%

Let $t\ge 2$ be an integer,
and $G$ be a~$2$-connected $K_{1,t}$-free graph.
We may assume that $t\ge 3$.
Let $F_0$ be a longest cycle in~$G$.
Since $G$ is $2$-connected,
we take a sequence of subgraphs $F_1, \dots, F_m$ such that 
each $F_\ell$ with $1\le \ell\le m$ is a maximum $2$-spider
among all $2$-spiders
from a vertex in~$G - \bigcup_{i=0}^{\ell-1} V(F_i)$
to $\bigcup_{i=0}^{\ell-1} V(F_i)$ and $V(G)=\bigcup_{i=0}^{m}V(F_i)$.
Then, $H:= \bigcup_{i=0}^m F_i$ is a~$2$-connected spanning subgraph in~$G$.

We show that $H$ is a desired $2$-connected spanning subgraph in~$G$.
Suppose that there is 
a vertex~$v \in V(H)$ with~$d_H(v)\ge t+1$.
Let $h_0$ be the smallest integer with~$v \in V(F_{h_0})$.
Note that $d_{F_{h_0}}(v) = 2$
and $d_{F_i}(v) \leq 1$ for~$i > h_0$.
Since $d_H(v) \geq t+1$,
there exist $t-1$ many $2$-spiders
$F_{h_1}, F_{h_2}, \ldots , F_{h_{t-1}}$
such that $h_i > h_0$ and $F_{h_i}$ contains $v$ as a foot for~$1 \leq i \leq t-1$.
Without loss of generality,
we may assume that $h_0 < h_1 < \cdots < h_{t-1}$.
Let~$x_0$ be either one of the two neighbors of~$v$ in $F_{h_0}$,
and~$x_i$ be the neighbor of~$v$ in~$F_{h_i}$ for~$1 \leq i \leq t-1$.
Since $G$ is $K_{1,t}$-free,
we see that $x_px_q \in E(G)$ for some $0 \leq p < q \leq t-1$.
Then $F_{h_p}-vx_p+x_px_q+vx_q$ is a~$2$-spider larger than $F_{h_p}$,
which contradicts the choice of~$F_{h_p}$.
Hence, we obtain the desired conclusion.
This completes the proof of Theorem~\ref{spider_mainthm}.
\qed

%%%%%%%%%%%%%%%%%%%%%%%%%%%%%%%%%%%%%%%%%%%%%%%%%%%%%%%%%%%%%%%
\subsection{Key lemma for the proof of Theorem~\ref{spider_mainthm}}
%%%%%%%%%%%%%%%%%%%%%%%%%%%%%%%%%%%%%%%%%%%%%%%%%%%%%%%%%%%%%%%

In the proof of Theorem~\ref{spider_mainthm},
we use an argument similar to the proof in the previous subsection,
but need to choose $j$-spiders more carefully.
To do that, we employ the following notation.

For a connected graph~$H$ and~$X \subseteq V(H)$,
if there is no cut set~$S$ with~$|S| \leq j-1$ in~$H$
such that $H-S$ has a component having no vertex in~$X$,
then we say that $H$ is {\it $j$-connected to~$X$}.
In other words,
$H$ is {\it $j$-connected to~$X$}
if and only if 
for any cut set~$S$ with~$|S| \leq j-1$ in~$H$,
each component of $H-S$ contains a vertex in~$X$.
Menger's theorem implies that 
if $|X| \geq j$ and a graph~$H$ is $j$-connected to~$X$,
then for any vertex $v \in V(H) \setminus X$,
there exists a~$j$-spider from~$v$ to~$X$.
We now show an important lemma for a graph $j$-connected to~$X$.

%%%%%%%%%%%%%%%%%%%%%%%%%%%%%%%%%%%%%%%%%%%%%%%%%%%%%%%%%%%%%%%
\begin{lemma}
\label{key}
Let $j \geq 2$ be an integer.
Let $H$ be a connected graph and~$X \subseteq V(H)$ with $V(H) \setminus X \neq \emptyset$
and $|X| \geq j$.
If $H$ is $j$-connected to~$X$, then
there is a vertex $w \in V(H)\setminus X$ such that $H-w$ is $j$-connected to~$X$
or
there is a~$j$-spider~$F$ from a vertex~$v$ in~$V(H)\setminus X$ to~$X$
such that $F$ contains all the neighbors in~$H$ of the branch $v$ of~$F$.
\end{lemma}
%%%%%%%%%%%%%%%%%%%%%%%%%%%%%%%%%%%%%%%%%%%%%%%%%%%%%%%%%%%%%%%
\pr
Let $j\ge 2$ be an integer,
$H$ be a connected graph and $X \subseteq V(H)$ with $V(H) \setminus X \neq \emptyset$ and $|X| \ge j$,
and assume that $H$ is $j$-connected to~$X$.

Suppose 
that
$H-w$ is not $j$-connected to~$X$ for any $w \in V(H)\setminus X$.
Since $H$ is $j$-connected to~$X$,
we see that
\begin{align}
\label{NOT}
\text{for any $w \in V(H) \setminus X$,
there exists a cut set $S$ in~$H$ with $|S| = j$ and $w \in S$}\\ \nonumber
\text{such that $H-S$ has a component having no vertex of~$X$.
} 
\end{align}
Since $V(H) \setminus X \neq \emptyset$,
there exists a cut set $S$ in~$H$ with~$|S| = j$,
and let $A$ be a component of~$H-S$ such that $A \cap X = \emptyset$.
(In this proof, we treat a component such as $A$ with its vertex set.)
We choose such $S$ and $A$ so that $|A|$ is as small as possible.
Let $\overline{A}=H-(S \cup A)$.
Note that $X\subseteq S \cup \overline{A}$.

We first claim that we may assume that $|A| \geq j+1$.
Suppose on the contrary that $|A| \leq j$.
Among all $j$-spiders from a vertex in $A$ to~$X$,
let 
$F$ be a maximum one.
Since $S$ is a cut set in~$H$ with $|S| = j$ and $A \cap X = \emptyset$,
the $j$-spider $F$ contains all the vertices in~$S$.
Let $u$ be the branch of~$F$.
We now show that $F$ contains all the neighbors of~$u$ in~$H$.
Suppose not, and let $v$ be a neighbor of $u$ in $H$ with $v \notin V(F)$.
In particular, $v \notin N_F(u)$.
Since $S \subseteq V(F)$, we see that $v \in A$.
If $N_H(v) \cap N_F(u) \neq \emptyset$, say $x \in N_H(v) \cap N_F(u)$,
then $F - ux + uv + vx$ is a $j$-spiders from a vertex in~$A$ to~$X$,
a contradiction to the maximality of~$F$.
Thus, 
we have $N_H(v) \cap N_F(u) = \emptyset$.
Since $H$ is $j$-connected to~$X$,
we have $|N_H(v)| \geq j$.
Since $u$ is the branch of the $j$-spider $F$,
we have $|N_F(u)| = j$.
These imply that 
$$2j \leq |N_H(v)| + |N_F(u)| = |N_H(v) \cup N_F(u)| \leq |(A \cup S)\setminus\{v\}|
= |A| + |S| -1 \leq 2j-1,
$$
a contradiction.
Therefore,
$F$ contains all the neighbors of~$u$ in~$H$,
and we are done.
This means that 
we may assume that $|A| \geq j+1$,
as claimed.

Let $u$ be a vertex in~$A$.
Since $u\notin X$,
it follows from (\ref{NOT}) that
there are a cut set $T$ in~$H$ with $|T|=j,~u \in T$ and $T\neq S$,
and a component $B$ of~$H-T$ such that $B\cap X=\emptyset$.
Let $\overline{B}=H-(T \cup B)$~(see Figure~1).
By the definitions of~$A$ and~$B$,
we see that $X\subseteq (S \cup \overline{A}) \cap (T \cup \overline{B})$.
Since $|A| \geq j+1$ and $|T|= j$,
we see that $A \cap B \neq \emptyset$ or $A \cap \overline{B} \neq \emptyset$.
We divide the proof into two cases.

%%%%%%%%%%%%%%%%%%%%%%%%%%%%%%%%%%%%%%%%%%%%%%%%%%
%%%%%%%%%%%%%%%%%%%%%%%%%%%%%%%%%%%%%%%%%%%%%%%%%%
% figure
%%%%%%%%%%%%%%%%%%%%%%%%%%%%%%%%%%%%%%%%%%%%%%%%%%
%%%%%%%%%%%%%%%%%%%%%%%%%%%%%%%%%%%%%%%%%%%%%%%%%%
\begin{figure}[h]
\center
\begin{tikzpicture}[samples=200,scale=0.8]

% Boxes
\draw (0,0)--(6,0)--(6,3)--(0,3)--cycle;
\draw (0,1)--(6,1);
\draw (0,2)--(6,2);
\draw (2,0)--(2,3);
\draw (4,0)--(4,3);

% Info
\draw (1,3.5-0.2)node{$A$};
\draw (3,3.5-0.2)node{$S$};
\draw (5,3.5-0.2)node{$\overline{A}$};
\draw (-0.5+0.1,0.5)node{$\overline{B}$};
\draw (-0.5+0.1,1.5)node{$T$};
\draw (-0.5+0.1,2.5)node{$B$};

\draw (1,0.5)node{$A\cap \overline{B}$};
\draw (1,1.5)node{$A\cap T$};
\draw (1,2.5)node{$A\cap B$};
\draw (3,0.5)node{$S\cap \overline{B}$};
\draw (3,1.5)node{$S\cap T$};
\draw (3,2.5)node{$S\cap B$};
\draw (5,0.5)node{$\overline{A}\cap \overline{B}$};
\draw (5,1.5)node{$\overline{A}\cap T$};
\draw (5,2.5)node{$\overline{A}\cap B$};

\end{tikzpicture}

%\label{fig}
\caption{Components of $H-S$ and $H-T$.}
\end{figure}
%%%%%%%%%%%%%%%%%%%%%%%%%%%%%%%%%%%%%%%%%%%%%%%%%%
%%%%%%%%%%%%%%%%%%%%%%%%%%%%%%%%%%%%%%%%%%%%%%%%%%

\noindent
{\bf Case 1:} $A \cap B \neq \emptyset$.

We set $S':=(T \cap (S \cup A)) \cup (S \cap B)$.
Then $A \cap B$ is a component of~$H-S'$ such that $(A\cap B)\cap X=\emptyset$ and
$|A \cap B|\le|A\setminus\{u\}|=|A|-1$ because $u\in T\cap A$.
If $|S'| \le j-1$, then
we get a contradiction to the assumption that $H$ is $j$-connected to~$X$;
if $|S'| = j$, then %since $w\in T \cap A$,
we get a contradiction to the choice of~$S$ and~$A$.
Hence $|S'|\ge j+1$.
We set $T':=(T \cap (S \cup \overline{A})) \cup (S \cap \overline{B})$.
Since $|S'|+|T'|=|S|+|T|=2j$,
we have
$|T'|= 2j-|S'|\le j-1$.
If $\overline{A} \cap \overline{B}\neq\emptyset$,
then $A\cup B$ is a component of~$H-T'$ such that $(A\cup B) \cap X=\emptyset$,
which contradicts the assumption that $H$ is $j$-connected to~$X$.
Thus $\overline{A} \cap \overline{B}=\emptyset$.
Then $X\subseteq T'$,
and hence %we have
$j\le |X|\le |T'|\le j-1$,
a contradiction.

\noindent
{\bf Case 2:} $A\cap B=\emptyset$ and $A \cap \overline{B} \neq \emptyset$.

We set $S'':=(T \cap (S \cup A)) \cup (S \cap \overline{B})$.
Then $A\cap\overline{B}$ is a component of~$H-S''$ such that $(A\cap \overline{B})\cap X=\emptyset$ and
$|A \cap \overline{B}|\le|A\setminus\{u\}|=|A|-1$ because $u\in T\cap A$.
If $|S''| \le j-1$, then
we get a contradiction to the assumption that $H$ is $j$-connected to~$X$;
if $|S''| = j$, then %since $w\in T \cap A$,
we get a contradiction to the choice of~$S$ and~$A$.
Hence $|S''|\ge j+1$.
We set $T'':=(T \cap (S \cup \overline{A})) \cup (S \cap B)$.
Since $|S''|+|T''|=|S|+|T|=2j$,
we have
$|T''|= 2j-|S''|\le j-1$.
If $\overline{A} \cap B\neq\emptyset$,
then $\overline{A}\cap B$ is a component of $H-T''$ such that $(\overline{A} \cap B)\cap X=\emptyset$,
which contradicts the assumption that $H$ is $j$-connected to~$X$.
Thus $\overline{A} \cap B=\emptyset$.
Since $A\cap B=\emptyset$,
we have $B\subseteq S$; thus $|B|\le|S|=j$.
However,
we have $|A| \geq j+1 > |B|$,
which contradicts the choice of~$S$ and~$A$.
This completes the proof of Lemma~\ref{key}.
\qed

%%%%%%%%%%%%%%%%%%%%%%%%%%%%%%%%%%%%%%%%%%%%%%%%%%%%%%%%%%%%%%%
%%%%%%%%%%%%%%%%%%%%%%%%%%%%%%%%%%%%%%%%%%%%%%%%%%%%%%%%%%%%%%%
\section{Proofs}
\subsection{Proof of Theorem~\ref{spider_mainthm}}
\label{proofspider_sec}
%%%%%%%%%%%%%%%%%%%%%%%%%%%%%%%%%%%%%%%%%%%%%%%%%%%%%%%%%%%%%%%
%%%%%%%%%%%%%%%%%%%%%%%%%%%%%%%%%%%%%%%%%%%%%%%%%%%%%%%%%%%%%%%

Let $j\ge 2$ and $t\ge2$ be integers,
and $G$ be a~$j$-connected $K_{1,t}$-free graph.
The case $j=2$ has been already proved in Subsection~\ref{proofj=2_sec},
so we may assume that $j \geq 3$.
Moreover, we may assume that $t\ge3$.

Let $F_0$ be a longest cycle in~$G$. 
Since $G$ is $j$-connected,
we have $|V(F_0)| \geq j$.
Let $G_0=G$ and $X_0 = V(F_0)$.
For $i \geq 1$, we define the graph $G_i$, the vertex $v_i \in V(G)$, 
the $j$-spider~$F_i$,
and the vertex subset $X_i$ of~$G$ recursively
so that $G_{i}$ is $j$-connected to~$X_{i}$
and $X_0 \subseteq X_1 \subseteq \cdots$, as follows:
Since $G$ is $j$-connected,
$G_0$ is $j$-connected to~$X_0$.
Suppose that $G_{i-1}, F_{i-1}$ and $X_{i-1}$ have been already defined,
and suppose further that $V(G_{i-1}) \setminus X_{i-1} \neq \emptyset$
and
$G_{i-1}$ is $j$-connected to~$X_{i-1}$.
\begin{enumerate}
\item[(A)]
Suppose that there is a vertex $v_i \in V(G_{i-1}) \setminus X_{i-1}$ such that
there is a~$j$-spider from~$v_i$ to~$X_{i-1}$ and it contains all the neighbors of~$v_i$ in~$G_{i-1}$.
Among such $j$-spiders, let $F_i$ be a maximum one.
Then,
let $G_i:=G_{i-1}$ and $X_i := X_{i-1} \cup V(F_i)$.

\item[(B)]
Otherwise, that is,
suppose that for all vertices~$v$ in~$V(G_{i-1})\setminus X_{i-1}$,
there is no $j$-spider from~$v$ to~$X_{i-1}$ such that it contains all the neighbors of~$v$ in~$G_{i-1}$.
Since $|X_{i-1}| \geq |X_{0}| \geq j$ and $G_{i-1}$ is $j$-connected to~$X_{i-1}$,
it follows from Lemma~\ref{key} that there exists a vertex~$v_i \in V(G_{i-1}) \setminus X_{i-1}$
such that $G_{i-1}-v_i$ is $j$-connected to $X_{i-1}$.
Since $G_{i-1}$ is $j$-connected to~$X_{i-1}$,
there exists a~$j$-spider from~$v_i$ to~$X_{i-1}$.
Among such $j$-spiders, let $F_i$ be a maximum one.
Then, let $G_i := G_{i-1} - v_i$ and $X_i := X_{i-1} \cup (V(F_i) \backslash \{v_i\})$.
\end{enumerate}
In either case,
since $X_{i-1} \subseteq X_i$,
it is easy to see that $G_{i}$ is $j$-connected to~$X_{i}$.
The above procedure continues until $X_m = V(G_m)$ holds
for some integer $m \geq 1$.
Note that $X_m \cup \{v_1, v_2, \dots , v_m\} = V(G)$, which is a disjoint union if (B) occurs in all $m$~steps.
We now show that 
every vertex~$v$ satisfies the desired properties.

First, 
suppose that 
there exists a~$j$-spider $F_{\ell}$ with~$1 \leq \ell \leq m$
whose branch is~$v$;
so, $v = v_{\ell}$.
Regardless of whether $v_\ell$ is chosen in~(A) or~(B),
it can be a foot of none of the $j$-spiders $F_0, F_1, \dots, F_m$,
since $N_{G_\ell}(v_\ell)\subseteq V(F_\ell)$ in case~(A), and $v_\ell\notin V(G_\ell)$ in case~(B).

Suppose next that 
there exists no $j$-spider $F_{\ell}$ with $1 \leq \ell \leq m$
whose branch is~$v$,
that is, $v \notin \{v_1,v_2,\ldots,v_m\}$,
and on the contrary that 
$v$ is a foot of $t-1$ of the $j$-spiders $F_0, F_1, \dots, F_m$,
say $F_{h_1}, F_{h_2}, \ldots, F_{h_{t-1}}$
with $h_1 < h_2 < \cdots < h_{t-1}$.
By the condition on~$v$,
there exists a~$j$-spider $F_{h_0}$ containing~$v$ as a non-branch vertex.
Note that $h_0 < h_i$ for $1 \leq i \leq t-1$.
Let $x_0$ be either one of the two neighbors of~$v$ in~$S_{h_0}$,
and $x_i$ be the neighbor of~$v$ in~$F_{h_i}$ for $1 \leq i \leq t-1$.
Since $G$ is $K_{1,t}$-free,
$x_px_q \in E(G)$ for some $0 \leq p < q \leq t-1$.
Then $F_{h_p}-vx_p+x_px_q+vx_q$ is a~$j$-spider from~$v_{h_p}$ to~$X_{h_p-1}$,
which contradicts the maximality of~$F_{h_p}$.
This completes the proof of Theorem~\ref{spider_mainthm}.
\qed

%%%%%%%%%%%%%%%%%%%%%%%%%%%%%%%%%%%%%%%%%%%%%%%%%%%%%%%%%%%%%%%
\subsection{Proof of Theorem~\ref{main}}
\label{proofmain_sec}
%%%%%%%%%%%%%%%%%%%%%%%%%%%%%%%%%%%%%%%%%%%%%%%%%%%%%%%%%%%%%%%

To prove Theorem~\ref{main},
we need a lemma concerning a factor theorem related to Hall's marriage theorem.
The following theorem is a consequence of known results;
it can be obtained by setting 
$n=j$,
$g(a) =0$,
$f(a) = \left\lceil \frac{1}{j}d_G(a) \right\rceil$ for each $a \in A$,
and $g(b) = f(b) = 1$ for each $b \in B$ 
in~\cite[Theorem 5.7]{AK-2011}.
This problem was originally raised
(with hints and a solution derived from a result of de Werra~\cite{deW}) in~\cite[Problem~7.11 on page~56]{L}.

%%%%%%%%%%%%%%%%%%%%%%%%%%%%%%%%%%%%%%%%%%%%%%%%%%%%%%%%%%%%%%%
\begin{theorem}[de Werra~\cite{deW}]%[Lov\'{a}sz~\cite{L}]
\label{factor}
Let $G=G(A,B)$ be a bipartite graph such that $d_G(b)=j$ for any $b \in B$.
Then there is a subgraph $H$ of~$G$
such that $d_{H}(a) \leq \left\lceil \dfrac{1}{j}d_G(a) \right\rceil$ for each $a \in A$
and $d_{H}(b)=1$ for each $b \in B$.
\end{theorem}
%%%%%%%%%%%%%%%%%%%%%%%%%%%%%%%%%%%%%%%%%%%%%%%%%%%%%%%%%%%%%%%

\noindent{\bf Proof of Theorem~\ref{main}.}
Let $j\ge 1$ and $k\ge j$ be integers,
and $G$ be a $j$-connected $K_{1,j(k-2)+2}$-free graph.
If $j=1$, then the theorem holds:
for $k=1$, $G$ consists of a single vertex;
for $k=2$, $G$ is complete;
and for $k\ge 3$, it follows from Theorem~\ref{JWthm2}.
Thus we may assume that $j\ge 2$.

By Theorem~\ref{spider_mainthm} with $t = j(k-2)+2$,
there exists a sequence of subgraphs $F_0, F_1, \dots, F_m$ in~$G$
such that $F_0$ is a cycle with $|V(F_0)| \geq j$,
each $F_{\ell}$ with $1 \leq \ell \leq m$
is a $j$-spider from a vertex in $G - \bigcup_{i=0}^{\ell-1} V(F_{i})$ to $\bigcup_{i=0}^{\ell-1} V(F_{i})$,
$V(G)=\bigcup_{i=0}^{m}V(F_i)$,
and every vertex $v$ satisfies the following property:
\begin{itemize}
\item
If $j \geq 3$ and there exists a~$j$-spider $F_{\ell}$~(with $1 \leq \ell \leq m$)
whose branch is~$v$,
then $v$ is a foot of none of the~$j$-spiders $F_0, F_1, \dots, F_m$.
\item
If $j=2$ or 
there exists no $j$-spider $F_{\ell}$~(with $1 \leq \ell \leq m$)
whose branch is~$v$, 
then $v$ is a foot of at most $j(k-2)$ of the $j$-spiders $F_0, F_1, \dots, F_m$.
\end{itemize}
Then, $H:= \bigcup_{i=0}^m F_i$ is a $2$-connected spanning subgraph in~$G$.
When $j \geq 3$, let
\[ W = \{v \mid \text{there exists a~$j$-spider $F_{\ell}$~(with $1 \leq \ell \leq m$)
whose branch is $v$}\}; \]
otherwise, that is, if $j=2$, let $W = \emptyset$.

Hereafter, we delete some edges from~$H$ to obtain a~$k$-tree of~$G$.
Let $K = K(A,B)$ be the bipartite graph
such that $A = V(H)$,
$B = \{F_{\ell} \mid 1 \leq {\ell} \leq m\}$,
and for $v \in A$ and $F_{\ell} \in B$, 
$vF_{\ell} \in E(K)$ if and only if
$F_{\ell}$ contains~$v$ as a foot. 
By the condition on $F_0, F_1, \dots, F_m$,
we have
\[ d_{K}(v) \leq
\begin{cases}
0 &\text{if $v \in W$,} \\
j(k-2) &\text{otherwise}
\end{cases}
\]
for $v\in A$.
For $F_{\ell} \in B$, 
since $F_{\ell}$ is a $j$-spider,
we have $d_{K}(F_{\ell}) = j$.
By applying Theorem~\ref{factor} to $K$,
we obtain a subgraph $L$ of $K$ such that
\[ d_{L}(v)
\leq \left\lceil \dfrac{1}{j}d_{K}(v) \right\rceil
\leq
\begin{cases}
0 &\text{if $v \in W$}, \\
k-2 &\text{otherwise}
\end{cases}
\]
for each $v \in A$, and $d_{L}(F_{\ell})=1$ for each $F_{\ell} \in B$.
Let $f$ be the function from $E(K)$ to~$E(H)$
such that for $vF_{\ell} \in E(K)$ with~$v \in A$ and $F_{\ell} \in B$,
$f(vF_{\ell})$ is the edge in~$H$ incident with~$v$ and contained in~$F_{\ell}$.
Let $T$ be the graph obtained from~$H$
by deleting all edges in~$f(E(K) \setminus E(\textcolor{red}{L}))$ and an edge in~$F_0$.
For each $1 \leq {\ell} \leq m$,
the vertices in~$V(F_{\ell}) \setminus \bigcup_{i=0}^{{\ell}-1} V(F_i)$ 
are connected to $\bigcup_{i=0}^{{\ell}-1} V(F_i)$ by the edge in $E(F_{\ell}) \cap f(E(L))$
because $d_{L}(F_{\ell})=1$,
and hence $T$ is connected.
We observe that $T$ has no cycle, and hence $T$ is a spanning tree of~$G$.
In addition, 
since $k \geq j$,
for each $v \in V(T)$,
\[ d_T(v)
= d_{F_{h_0}}(v) + d_{\textcolor{red}{L}}(v)
\leq
\begin{cases}
j+ d_{\textcolor{red}{L}}(v) \leq j+0 \leq k &\text{if $v \in W$,} \\
2 + d_{\textcolor{red}{L}}(v) \leq k &\text{otherwise,}
\end{cases}
\]
where $h_0$ is the smallest integer with $v \in V(F_{h_0})$.
This completes the proof of Theorem~\ref{main}.
\qed

%%%%%%%%%%%%%%%%%%%%%%%%%%%%%%%%%%%%%%%%%%%%%%%%%%%%%%%%%%%%%%%
\subsection{Arguments related to Theorem \ref{maxdeg_thm}}
%%%%%%%%%%%%%%%%%%%%%%%%%%%%%%%%%%%%%%%%%%%%%%%%%%%%%%%%%%%%%%%

As we have explained in Subsection~\ref{introtree_sec},
an weaker version of Theorem~\ref{maxdeg_thm} can be proved
by the same arguments as in the previous subsection.
Note that in the sequence of subgraphs $F_0, F_1, \dots, F_m$
obtained from Theorem~\ref{spider_mainthm},
every vertex~$v$ satisfies the following property:
\begin{itemize}
\item
If $j \geq 3$ and there exists a~$j$-spider $F_{\ell}$~(with $1 \leq \ell \leq m$)
whose branch is $v$,
then $v$ is a foot of none of the $j$-spiders $F_0, F_1, \dots, F_m$.
\item
If $j=2$ or 
there exists no $j$-spider $F_{\ell}$~(with $1 \leq \ell \leq m$)
whose branch is~$v$, 
then $v$ is a foot of at most $d_G(v) - 2$ of the $j$-spiders $F_0, F_1, \dots, F_m$.
\end{itemize}
Then,
following the proof of Theorem~\ref{main},
we obtain a spanning tree~$T$
such that $d_T(v) \leq \max\left\{\left\lceil \frac{d_G(v) -2}{j}\right\rceil +2, j\right\}$
for any vertex~$v$.
We leave the detail for the readers.

%%%%%%%%%%%%%%%%%%%%%%%%%%%%%%%%%%%%%%%%%%%%%%%%%%%%%%%%%%%%%%%
%%%%%%%%%%%%%%%%%%%%%%%%%%%%%%%%%%%%%%%%%%%%%%%%%%%%%%%%%%%%%%%
\section*{Declaration of competing interest}
%%%%%%%%%%%%%%%%%%%%%%%%%%%%%%%%%%%%%%%%%%%%%%%%%%%%%%%%%%%%%%%
%%%%%%%%%%%%%%%%%%%%%%%%%%%%%%%%%%%%%%%%%%%%%%%%%%%%%%%%%%%%%%%

The authors declare that they have no known competing financial interests or personal relationships that could have appeared to influence the work reported in this paper.

%%%%%%%%%%%%%%%%%%%%%%%%%%%%%%%%%%%%%%%%%%%%%%%%%%%%%%%%%%%%%%%
%%%%%%%%%%%%%%%%%%%%%%%%%%%%%%%%%%%%%%%%%%%%%%%%%%%%%%%%%%%%%%%
\section*{Data availability}
%%%%%%%%%%%%%%%%%%%%%%%%%%%%%%%%%%%%%%%%%%%%%%%%%%%%%%%%%%%%%%%
%%%%%%%%%%%%%%%%%%%%%%%%%%%%%%%%%%%%%%%%%%%%%%%%%%%%%%%%%%%%%%%

No data was used for the research described in the article.

%%%%%%%%%%%%%%%%%%%%%%%%%%%%%%%%%%%%%%%%%%%%%%%%%%%%%%%%%%%%%%%
%%%%%%%%%%%%%%%%%%%%%%%%%%%%%%%%%%%%%%%%%%%%%%%%%%%%%%%%%%%%%%%
\section*{Acknowledgments}
%%%%%%%%%%%%%%%%%%%%%%%%%%%%%%%%%%%%%%%%%%%%%%%%%%%%%%%%%%%%%%%
%%%%%%%%%%%%%%%%%%%%%%%%%%%%%%%%%%%%%%%%%%%%%%%%%%%%%%%%%%%%%%%

This work was supported by the Research Institute for Mathematical Sciences, an International Joint Usage/Research Center located in Kyoto University.
The authors would like to thank Prof.~Kano
for informing them about the reference for Theorem~\ref{factor}.

%%%%%%%%%%%%%%%%%%%%%%%%%%%%%%%%%%%%%%%%%%%%%%%%%%%%%%%%%%%%%%%
%%%%%%%%%%%%%%%%%%%%%%%%%%%%%%%%%%%%%%%%%%%%%%%%%%%%%%%%%%%%%%%

\end{document}